\documentclass[]{siamart220329}
\usepackage{amsmath}
\usepackage{amssymb}
\usepackage{graphicx}
\usepackage{tikz}
\usetikzlibrary{external}
\usepackage{pgfplots}
\usepackage[Publish]{externaltikz}
\usepackage{calc}

\newcommand{\R}{\mathbb{R}} 
\pgfplotsset{compat=1.18}
\newtheorem{remark}{Remark}
\begin{document}

\title{A spectral approach to interface  layers on networks for the linearized BGK equation and its acoustic limit}

\author{R. Borsche\footnotemark[1] 
	\and T. Damm \footnotemark[1] 
     \and  A. Klar\footnotemark[1] 
 \and Y. Zhou\footnotemark[2] }
\footnotetext[1]{RPTU  Kaiserslautern, Department of Mathematics, 67663 Kaiserslautern, Germany 
(\{borsche, damm, klar\}@rptu.de)}
\footnotetext[2]{RWTH Aachen, Department of Mathematics, Aachen, Germany 	(zhou@igpm.rwth-aachen.de)}

 
\date{}


\maketitle

\begin{abstract}
	We consider in this paper a velocity discretized version of the full   linear kinetic BGK model and the corresponding limit for small Knudsen number, the linearised Euler or acoustic system. Considering these equations 
	on networks,
	coupling conditions for the macroscopic equations are derived from the kinetic conditions via an asymptotic analysis near the nodes of the network. Here, a  degeneracy in the limit equations requires not only the investigation of kinetic layers, but also the discussion of viscous
	 layers.  Using the kinetic coupling conditions at the junction and coupling  kinetic  and viscous layers   to the outer problems on the edges
	 one obtains 
a  coupled kinetic half-space problem at each node. 
	A spectral  method is  developed to solve this coupled  kinetic half-space problems. 
	This allows to obtain a detailed picture of the various interface layers near the nodes and to determine the relevant coefficients in the 
	kinetic derived coupling conditions for the macroscopic equations.
	Numerical results show the accuracy and efficiency  of the approach. 
\end{abstract}

{\bf Keywords.} 
Kinetic layer, coupling condition, kinetic half-space problem, networks

{\bf AMS Classification.}  
82B40, 90B10,65M08


\section{Introduction}

There have been many attempts  to define coupling conditions for macroscopic partial differential equations on networks including, for example,  drift-diffusion equations, scalar  hyperbolic equations, or hyperbolic systems like  the wave equation or Euler type models, see for example \cite{BGKS14,CC17, BNR14,BHK06a,BHK06b,CHS08,EK16,BCG10,HKP07,CGP05,G10}. 
On the other hand, coupling conditions for kinetic equations on networks have been discussed in a much smaller number of publications,
see \cite{FT15,HM09,holle,BKKP16}. 
In \cite{BKKP16} a first attempt to derive a coupling condition for a
macroscopic equation from the underlying kinetic model has been presented for the case of a kinetic equations for chemotaxis.
In \cite{BK18,ABEK22} more  general and more accurate approximation procedures have been presented and discussed for linear kinetic equations.
They are motivated by the classical procedure to find kinetic slip boundary conditions for macroscopic equations based on the analysis of the 
kinetic layer \cite{BSS84,BLP79,G92,G08,UTY03}.
In a recent work \cite{BDKZ1} the non-degenerate case has been reconsidered, existence and uniqueness of the kinetic coupled half-pace problems has been proven and a general numerical approach has been presented for the velocity discretized version of the kinetic equations.

In the present  work, we consider the degenerate case of the full linear BGK equation. As in \cite{BDKZ1} we consider 
the velocity-discretized equations and 
derive coupling conditions for the associated macroscopic equations 
on a network from the  kinetic conditions via an asymptotic analysis of the situation near the nodes. 
The main additional difficulties compared to \cite{BDKZ1}  are due to the degeneracy of the limit equations. This requires the additional consideration of the viscous layers on the edges and their coupling to the bulk solution on the one hand and the kinetic layer on the other hand.
We concentrate here on the modelling and numerical aspects of the problem. A thorough analytical investigation, including a formal asymptotic theory to  higher orders and a convergence result will be given in \cite{ZK}.

The paper is organized in the following way. 
In section \ref{sec:equations} we discuss  the kinetic and macroscopic equations on networks considered here and  classes of coupling conditions for these equations.
In Section \ref{kinlayer}     an 
asymptotic analysis of the kinetic equations near the nodes and resulting  kinetic layers at the nodes are discussed as well as 
 the associated viscous layer.
This leads to an abstract formulation of the coupling conditions for the macroscopic equations at the nodes based on  coupled kinetic half-space problems. 
In the following Section \ref{discretekin} a general velocity discretization of the kinetic equation via kinetic discrete velocity models is considered and the associated  kinetic moment problem is given.
Moreover, the discrete layer problem in moment coordinates is investigated and a spectral method is discussed.
In Section \ref{couplingconditions} the coupling conditions for the macroscopic equations as well as details on the kinetic solution near the node in the 
asymptotic limit  are obtained.
Moreover, the section includes a short overview of the results of two approximate methods to determine the coupling conditions.
Finally, Section 
\ref{numericalbgk} illustrates the above findings for 
 a simple tripod network comparing kinetic and macroscopic network solutions and giving  numerical results for various examples.

\section{Equations and coupling conditions}
\label{sec:equations}

For  $f_\epsilon=f_\epsilon(x,v,t)$ with $x\in\R$ and  $v\in\R$ at time $t\in[0,T]$ we consider the linear kinetic BGK model with a hyperbolic space-time scaling 
\begin{equation}\label{eq:2.1bgk}
\partial_t f_\epsilon +v\partial_xf_\epsilon=\frac{1}{\epsilon}Q(f_\epsilon)= -\frac{1}{\epsilon}\left(f_\epsilon-\left(\rho+vq+\frac12(v^2-1)(S-\rho)\right)M(v)\right),
\end{equation} 
where density, mean flux and total energy are given by 
$$
\rho_\epsilon=\int_{-\infty}^{\infty}f_\epsilon(v)dv, \quad q_\epsilon=\int_{-\infty}^{\infty}vf_\epsilon(v)dv, \quad S_\epsilon=\int_{-\infty}^{\infty}v^2f_\epsilon(v)dv
$$
and the standard Maxwellian is defined by
$$
M(v)=\frac{1}{\sqrt{2\pi}}\exp(-\frac{v^2}{2}).
$$
Determining the equations for  the first 4 moments of the kinetic  system gives 
\begin{align}
\begin{aligned}
\partial_t  \rho_\epsilon+    \partial_x  q_\epsilon&= 0\\
\partial_t  q_\epsilon+   \partial_x  S_\epsilon &= 0\\
\partial_t  S_\epsilon+     \partial_x   h_\epsilon &= 0\\
\epsilon \partial_t h_\epsilon +  \epsilon \partial_x \int v^4 f_\epsilon dv   &= - ( h_\epsilon -3 q_\epsilon).
\end{aligned}
\end{align}
where 
$ h_\epsilon = \int v^3  f_\epsilon dv $.
To order $\epsilon$ one obtains 
\begin{align}
\begin{aligned}
h_\epsilon= 3 q_\epsilon+ \mathcal{O}(\epsilon)\\
\end{aligned}
\end{align}
and thus, the  associated limit equation for $\epsilon\to 0$ as
\begin{equation}
\begin{aligned}
\partial_t\rho_0+\partial_xq_0&=0, \\
\partial_tq_0+\partial_xS_0&=0,\\
\partial_tS_0 +a^2\partial_xq_0&=0 
\end{aligned}
\label{eq:2.2}
\end{equation}
with $a^2=3$. The equations are called the acoustic system.
Moreover,
\begin{equation}
	\begin{aligned}
f_0 =  \left( \rho_0+v q_0+\frac12(v^2-1)( S_0-\rho_0)\right)M(v) .
\end{aligned}
\end{equation}
The characteristic variables for equations  \eqref{eq:2.2} associated to the eigenvalues $\lambda_\mp= \mp a$ and
$\lambda_0 =0$ are 
\begin{eqnarray}\label{eq:nribgk}
\begin{aligned}
r_{-}&=& S_0-aq_0,\qquad 
r_0&=& S_0-a^2\rho_0,\qquad 
r_{+}&=&S_0+aq_0.
\end{aligned}
\end{eqnarray}
See \cite{BGL00,Bel04,GL01} for the derivation of the acoustic limit from  Boltzmann and  BGK equation.
We note that (\ref{eq:2.2}) is the linearisation of the compressible Euler equations around a state with velocity zero.



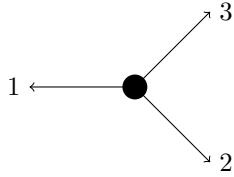
\begin{figure}[h!]
	\center
	\externaltikz{NodeSketch}{
		\begin{tikzpicture}
		\draw[->] (0,0)--(-1.4,0) node[left]{$1$};
		\draw[->] (0,0)--(1,-1) node[right]{$2$};
		\draw[->] (0,0)--(1,1) node[right]{$3$};
		\node[fill=black,circle] at (0,0){};
		\end{tikzpicture} 
}
	\caption{Node connecting three edges and orientation of the edges.}
	\label{fig:tripod}
\end{figure}
If the kinetic equation is considered on a network, it is sufficient to study a single node, see Figure \ref{fig:tripod}.
At each node so called coupling conditions are required. 
In the following we consider a node connecting $n$ edges, which are oriented away from the node, as in Figure \ref{fig:tripod}.
Each edge $i$ is parametrized by the interval $[0,b_i]$ and the kinetic and macroscopic quantities are denoted by
$f_\epsilon^i$ and $\rho_0^i, q_0^i,S_0^i$ respectively.
On the kinetic level for each edge a condition on  the ingoing part  $f^i(0,v),v>0$ of the distribution function  is required at the node, i.e.
at   $x=0$.
For the network problem  a possible choice of such a coupling condition for the kinetic problem is
\begin{align}
\label{kincoup}
f_\epsilon^i(0,v) = \sum_{j=1}^n \beta_{ij}f^j_\epsilon (0,-v), v >0\ , i= 1, \ldots ,n, 
\end{align}
compare \cite{BKKP16}.
Then, the total mass in the system is conserved, if 
\begin{align}\label{eq:conservative_coupling_f}
\sum_{i=1}^n \beta_{ij} = 1,
\end{align}
since in this case the balance of fluxes,  i.e. $\sum_{j=1}^n \int v f^j_\epsilon(0,v) dv =0$, holds.
\begin{remark}
Additionally, the coupling conditions also yield for any $k =2l+1$ with $l= 1 , \ldots$
$$\sum_{j=1}^n \int v^k f^j_\epsilon(0,v) dv =0  $$
\end{remark}

In particular, we will consider the case of a node with symmetric coupling conditions, 
that means  $\beta_{ij} = \frac{1}{n-1}, i \neq 0$  and $\beta_{ii} =0$.

The  coupling conditions for the macroscopic equations  on the network are  conditions for the $3n$ macroscopic quantities $(\rho_0^i(0),q_0^i(0),S_0^i(0))$ at the nodes $x=0$.
They   are given by $n$ 
conditions
for  the ingoing characteristic variables at the nodes 
\begin{align*}
r_+^i (0) = S_0^i(0) + a q_0^i(0).
\end{align*}
These $n$ conditions are given by the balance of fluxes
$$\sum_{i=1}^n q_0^i (0) =0$$
and, in particular in case of symmetric nodes, 
 by invariants at the nodes. For the present system of 3 equations with only one ingoing characteristic 
we need  only one more invariant at the node, for example   $S_0^i(0) + \delta q_0^i (0)$ with an arbitrary $\delta  >0$  giving  $n-1$ additional
equations
$$
S_0^i(0) + \delta q_0^i (0)=S_0^j(0) + \delta q_0^j(0)
$$
for $i,j=1, \ldots,n$.
This is combined with $n$ further conditions given by the  prescribed values of the outgoing characteristics at the nodes
\begin{align}
	\label{outcond}
S_0^i(0) - a q_0^i(0) = r^i_- (0) .
\end{align}

The conditions discussed up to now give a sufficient number of conditions for the macroscopic equations, since there is
no condition required for the 
degenerate characteristic $r_0$ associated to the eigenvalue $\lambda_0=0$.

\begin{remark}
	We mention, that the number  $\delta$ plays a similar role as the so-called extrapolation length for kinetic boundary layers, \cite{BSS84}.
\end{remark}

In the following we will not only determine the coupling conditions for the macroscopic equations at the node, but also give a complete analysis of the kinetic solution and the structure of the interface layers near the node  in the limit $\epsilon \rightarrow 0$.

\section{Kinetic and viscous layers  at the nodes and coupling conditions for macroscopic equations}
\label{kinlayer}


\subsection{The kinetic layer}
We perform a  kinetic layer analysis, see \cite{BSS84,C69,CGS88,Bardos,GMP}, at the left boundary of  each edge $[0,b_i]$.
A rescaling of the spatial variable in equation (\ref{eq:2.1bgk}) with $\epsilon$ results  in  the scaled equation
$$
\partial_t f +  \frac{1}{\epsilon} v \partial_x f = \frac{1}{\epsilon} Q(f).
$$
This yields to zeroth order in $\epsilon$ to the  following stationary kinetic  half space  problem for 
the scaled spatial variable $x \in [0,\infty]$
\begin{align}\label{bgkhalfspacebgk}
\begin{aligned}
v \partial_x f =\frac{1}{\epsilon}Q(f)= -\frac{1}{\epsilon}\left(f-\left(\rho+vq+\frac12(v^2-1)(S-\rho)\right)M(v)\right), 
\end{aligned}
\end{align}
where $\rho,q$ and $S$ are again the zeroth, first and second moments of  $f$. 
At $x=0$ one has to prescribe for such  half space problem, as for the original kinetic problem,  the ingoing characteristics, i.e.
$$
  f(0,v), v > 0\ .
$$
At  $x=\infty$,  a single further condition is   needed. It is in the present case given by the outgoing macroscopic characteristics as 
\begin{align}
\label{infcond}
 \int \left(v^2 - a v \right)fd v= r_-(0)= S_0(0)-a q_0(0).
\end{align}
%
The resulting solution of the half -space problem at infinity has the form 
$$
f(\infty,v) =\left(\rho_\infty +vq_\infty+\frac12(v^2-1)(S_\infty-\rho_\infty)\right)M(v)
$$

%

The resulting outgoing  solution of the half space problem at $x=0$ is 
$
f(0,v), v<0\ 
$
or
$
f(0,-v), v > 0\ .
$
We combine now the layers on all edges adjacent to the node under consideration and use the kinetic coupling conditions
\begin{align*}
f^i(0,v) = \sum_{j=1}^n \beta_{ij}f^j (0,-v), v >0.
\end{align*}
We use the kinetic  coupling conditions to obtain conditions on the ingoing solutions of the half space problems on the different arcs.
That means 
\begin{align*}
f^i(0,v)
= \sum_{j=1}^n \beta_{ij} f^j(0,-v) , v>0
\end{align*}
Thus, we obtain  a coupled kinetic half-space problem at each  node. 
The resulting values at infinity determined from the coupled half-space problems finally give the desired coupling conditions
for the macroscopic quantities.

However, it is easy to note that this problem does not have a unique solution in the present situation.
Any linear combination of collision invariants $(\alpha + \beta v  + \gamma v^2)M$ solves the linear equations.
Moreover, using  the same values of $\alpha, \beta, \gamma $  for all edges it fulfils also the coupling 
conditions, if $ \beta =0$. The limit conditions 
$$
\int \left(v^2 - a v \right)fd v= 0
$$
are fulfilled for all edges, if 
$$
\int \left(v^2 - a v \right ) (\alpha +  \gamma v^2) M d v= 0,
$$
i.e. if $\alpha = - 3 \gamma.$
In other words  for any solution $f=(f^1, \cdots f^n)$ of the coupled half-space problem the function $\tilde f= \varphi+C (v^2-3 ) M(v)$  
for an arbitrary constant $C\in \R$ is also a solution fulfilling (\ref{infcond}).
To fix this constant  we have to prescribe an additional condition.
This is done by considering the full situation at the node including the additional viscous layer.

\subsection{The viscous layer}
\label{visclayer}
In this subsection we consider the viscous layer in more detail.
On each edge a viscous layer of size of order $\sqrt{\epsilon}$ connects  the kinetic layer solution to the bulk solution on the edge. 
The associated  scaling for the viscous  layer is $x \rightarrow \frac{x}{\sqrt{\epsilon}}$.
We write in the viscous layer
$$f (x,v,t)= f_0  (x,v,t) + \hat f \left(\frac{x}{\sqrt{\epsilon}},v,t\right),$$ where $\hat f$ fulfils the scaled equation
$$
\partial_t \hat f +  \frac{1}{\sqrt{\epsilon}} v \partial_x \hat f = \frac{1}{\epsilon} Q(\hat f).
$$
This is equivalent to the moment system 
\begin{align}
\begin{aligned}
\partial_t \hat \rho+  \frac{1}{\sqrt{\epsilon}}  \partial_x \hat q&= 0\\
\partial_t \hat q+   \frac{1}{\sqrt{\epsilon}}  \partial_x \hat S &= 0\\
\partial_t \hat S+   \frac{1}{\sqrt{\epsilon}}    \partial_x  \hat h &= 0\\
\epsilon \partial_t\hat h +  \sqrt{\epsilon} \partial_x \int v^4 \hat f dv   &= - (\hat h -3 \hat q),
\end{aligned}
\end{align}
where $\hat \rho, \hat q, \hat S$ are the corresponding moments of $\hat f$ and $\hat h = \int v^3 \hat f dv $.
To order $\sqrt{\epsilon}$ one obtains 
\begin{align}
\begin{aligned}
\sqrt{\epsilon}\partial_t \hat \rho+  \partial_x \hat q&= 0\\
\sqrt{\epsilon}\partial_t \hat q+    \partial_x \hat S &= 0\\
\sqrt{\epsilon}\partial_t \hat S  + 3 \partial_x \hat q &=   \sqrt{\epsilon} \partial_{xx} \int v^4 \hat f dv .
\end{aligned}
\end{align}
Using $\hat f= \left(\hat \rho+v\hat q+\frac12(v^2-1)(\hat S-\hat \rho)\right)M(v) + \mathcal{O} (\sqrt{\epsilon})$ we obtain
again to order $\sqrt{\epsilon}$ for the last equation
\begin{align}
\begin{aligned}
\sqrt{\epsilon}\partial_t \hat S +3 \partial_x \hat q  &=  \sqrt{\epsilon} \partial_{xx} \int v^4   \left(\hat \rho+v\hat q+\frac12(v^2-1)(\hat S-\hat \rho)\right)M(v)dv 
\end{aligned}
\end{align}
and
\begin{align}
\begin{aligned}
\sqrt{\epsilon}\partial_t \hat S  +3 \partial_x \hat q &=  \sqrt{\epsilon} \partial_{xx}   \left(6\hat S-3 \hat \rho\right)
\end{aligned}
\end{align}
and then
\begin{align}
\begin{aligned}
\partial_t (\hat S   - 3 \hat \rho)&=  \partial_{xx}   \left(6\hat S-3 \hat \rho\right).
\end{aligned}
\end{align}
From $\sqrt{\epsilon}\partial_t \hat q+    \partial_x \hat S =0$ we have  for $\epsilon \rightarrow 0$ in the viscous layer $\partial_x \hat S= \mathcal{O} (\sqrt{\epsilon})$ and therefore to $ \mathcal{O} (\sqrt{\epsilon})$
\begin{align}
\begin{aligned}
\partial_t (\hat S   - 3  \hat \rho)&=  \partial_{xx}   \left( \hat S-3\hat  \rho\right)
\end{aligned}
\end{align}
or $\partial_t \hat r =  \partial_{xx}  \hat r $ with $\hat r = \hat S   - 3  \hat \rho$.

This is a diffusion equation for the characteristic associated to the zero eigenvalue in the macroscopic equations.
The third order moment in the viscous layer is approximated by
\begin{align}
\begin{aligned}
\hat h = \int v^3 \hat f dv= 3 \hat q -\sqrt{\epsilon} \partial_x( \hat S - 3 \hat \rho) + \mathcal{O} (\epsilon)
= 3 \hat q -\sqrt{\epsilon} \partial_x \hat r + \mathcal{O} (\epsilon).
\end{aligned}
\end{align}
Matching  the viscous layer  to the bulk solution yields  the boundary condition at infinity for the viscous layer solution as 
 \begin{align*}
  \hat r(\infty,t) = 0 .
 \end{align*}
Matching the viscous layer equations to the kinetic layer gives the following condition
\begin{align*}
S(\infty,t)-3 \rho(\infty,t) = r_0(0,t) + \hat r  (0,t) 
\end{align*}
where $S(\infty,t)=S_\infty (t), \rho(\infty,t)= \rho_\infty(t)$ are the asymptotic states of the kinetic layer

%

\subsection{The final coupling condition}
\label{vfinalcoup}
In this subsection we  obtain the missing condition for the kinetic half-space problems.
As mentioned before,  the kinetic coupling conditions yield 
$$
\sum_{i=1}^n\int v^3 f^i (0,v,t) dv  =0,
$$
where $f^i(x,v,t)$ denotes the distribution function on each edge.
Using the invariance of $\int v^3 f(x,v,t) d v $ in the kinetic layer, 
we obtain the following  $\mathcal{O} ({\epsilon})$ approximation of the kinetic third order moment at $x=0$
\begin{align*}
\int v^3 f(0,v,t) d v &=\int v^3 f(\infty ,v,t) d v  =\int v^3 f_0(0 ,v,t) d v + \int v^3 \hat f(0 ,v,t) d v \\=& h_0(0,v,t)+\hat h (0,v,t) =
 3 q_0 + 3 \hat q - \sqrt{\epsilon} \partial_x  \hat r(0,t)+ \mathcal{O}(\epsilon).
\end{align*}
Using now the the balance of fluxes  and the invariance of $\int v f(x,v,t) d v$ in the kinetic layer we obtain
$$
0= \sum_{i=1}^n\int v f^i (0,v,t) dv  = \sum_{i=1}^n\int v ( f_0^i (0,v,t) +  \hat f^i (0,v,t) ) dv    =\sum_{i=1}^n (q_0^i + \hat q^i ) =  0.
$$
This leads to the condition
$$
\sum_{i=1}^n \partial_x \hat r^i(0,t)  =0.
$$
Considering now the viscous layer solution for the quantity
\begin{align*}
\hat R(x,t) = \sum_{i=1}^n  \hat r^i(x,t) 
\end{align*}
we obtain 
\begin{align}
\begin{aligned}
\partial_t \hat R&=  \partial_{xx}  \hat R
\end{aligned}
\end{align}
 with the conditions $\hat R(\infty,t) =0$ and $\partial_x \hat R(0,t) =0$.
 This gives the trivial solution $\hat R(x,t) =0$ and the desired condition
\begin{align}
\label{balancerho}
\sum_{i=1}^n (S^i(\infty,t)- 3 \rho^i(\infty,t) ) = R_0(0,t)+ \hat R(0,t) = R_0(0,t)= \sum_{i=1}^n  (S^i_0(0,t) - 3   \rho^i_0(0,t))
\end{align}
for the coupled kinetic half-space problem.
In Section \ref{couplingconditions} the solution of  coupling problem using the above discussed coupling condition is considered  in detail for the case of a kinetic discrete velocity model
with arbitrary numbers of velocities.

%


\subsection{The situation at the node}
\label{nnode}

Considering first the variables $q$ and $S $, the quantities near the node  are given by the solution $q_\infty$ and $S_\infty$ at infinity  of the kinetic layer problems. They are both constant throughout the kinetic layer, since
in the layer $\int v \varphi dv $ as well as $\int v^2 \varphi dv $ are constant.
Moreover, $q$ and $S$ are also constant throughout the viscous layer due to the considerations in section \ref{visclayer}.
Finally, a wave travelling  with speed $a$ is observed in $S$ and $q$ with left states  $q_\infty$ and $S_\infty$ and right states given by the
corresponding  initial values.

The situation is more complex for the variable $\rho$. 
For  this  quantity  kinetic and viscous  layers  may  appear  additionally to the travelling wave solution.
Since there is no travelling wave in the zero characteristic $S-3\rho$ and since the left state of the travelling wave solution in $S$ 
is $S_\infty(t)$, we obtain the  left state of the travelling wave solution in $\rho$ 
as 
$$
 \rho_0 (0,t) + \frac{1}{3}(S_\infty (t)- S_0(0,t)) .
$$
Consequently, the solution in $\rho$ consists of a kinetic layer from $\rho(0,t)$, the solution at the node to $\rho_\infty(t)=\rho(\infty,t)$
at the end of the kinetic layer. Then, there is a  viscous layer near the node ranging from $\rho_\infty (t)$  to  the left state of the travelling wave solution
computed above.

\section{The discrete velocity BGK equation}
\label{discretekin}

For a more detailed analysis of the coupled layer problem at the node and a determination of the macroscopic 
coupling coefficient $\delta$ we concentrate on a discrete velocity version of the kinetic equations, see \cite{Bab,Ber}
for more details on discrete boundary layers.
We proceed as follows.
As in the work of F. Coron \cite{coron} we consider  orthonormal Hermite  polynomials $P_k(v), k=0, \ldots ,2N$ on $[-\infty,\infty]$ defined by $P_0=\frac{1}{\pi^{1/4}},P_1 = \frac{\sqrt{2}}{\pi^{1/4}}v$ and
\begin{align*}
	v P_k(v) = \alpha_{k+1} P_{k+1} + \alpha_k P_{k-1}, k=1, \cdots 2N-1
\end{align*}
with $\alpha_k = \sqrt{\frac{k}{2}}  $, compare again \cite{coron}.
Note that $P_2= \sqrt{2} v^2P_0 -\frac{1}{\sqrt{2}}P_0$, $P_3= \frac{2}{\sqrt{3}} v^3P_0 - \frac{\sqrt{3}}{\sqrt{2}}P_1$, $  P_4=  \sqrt{\frac{2}{3}}v^4 P_0 - \sqrt{3} P_2 - \frac{\sqrt{3}}{2 \sqrt{2}}P_0 $.
Moreover, Define the associated functions $$H_k = P_k \exp(-\frac{v^2}{2}).$$
Using  the  transformations $v = \sqrt{2} \tilde v $  and $f= \tilde f H_0$
the kinetic equation can be rewritten as 
\begin{equation}
	\partial_t f +\sqrt{2} v\partial_xf= -\frac{1}{\epsilon}\left(f-\left(H_0 g_0 +H_1 g_1 + H_2 g_2\right)\right),
\end{equation} 
with 
$$
g_k= \int H_k(v)  f(v) dv , k=0,1,2.
$$
Kinetic density, mean flux and total energy $\rho,q,S$ are then obtained from 
$$
\rho=\sqrt{2}g_0, \quad q= \sqrt{2} g_1 \quad S=2 g_2+\rho
$$

\subsection{The discrete velocity model}
\label{dvmbgk}
For the following we discretize the BGK-equation (\ref{eq:2.1bgk}) in velocity space and obtain a kinetic discrete velocity model  for the discrete distribution functions 
$f_i(x,t), i= 1, \cdots  2N$ as 

\begin{align}\label{eq:kineticbgk}
	\begin{aligned}
		\partial_t f_i + v_i\partial_xf_i &= -\frac{1}{\epsilon}\left(f_i-M_i\right)
	\end{aligned}
\end{align}
with  the velocity discretization $$-1\le v_1< v_2 < \cdots < v_N < 0 < v_{N+1}< \cdots < v_{2N-1} < v_{2N} \le 1.$$
We  assume for  symmetry 
\begin{align*}
	v_{2N} = -v_1, 
	\cdots, 
	v_{N+1}= -v_N.
\end{align*}
Let $w_i \ge 0, i=1, \ldots 2N$ be symmetric weights, such that $\ \sum_{i=1}^{2N}  w_i =1$.
We choose $v_i, i=1, \ldots 2N$ to be the Gauß-Hermite  points on $[-\infty,\infty]$ and $w_i$ the associated 
Gauss-Hermite weights.

Moreover, we use a  discrete  linearised Maxwellian $M_i$ given by
\begin{align}
	M_i = w_i e^{v_i^2} ( H_0 (v_i)g_0 + H_1(v_i ) g_1 + H_2(v_i)g_2)
\end{align}
with $g_k,k=0, \ldots, 2N-1$ defined by
\begin{align*}
	g_k  = \sum_{i=1}^{2N} H_k(v_i) f_i\ .
\end{align*}
The choice of discrete Maxwellian yields for $k=0,1,2$
\begin{align}
\label{mommaxbgk}
\sum_{i=1}^{2N}M_i H_k(v_i) = g_k
\end{align}
and
\begin{align}
\sum_{i=1}^{2N}M_i H_k(v_i) = 0, k=3, \ldots,2N-1
\end{align}
due to discrete orthogonality.

Moreover,  
\begin{align*}
\sum_{i=1}^{2N}H_{2N} (v_i)f_i\ =0.
\end{align*}

Let now $G=(u,g) T$ with $u=(g_0,g_1,g_2,g_3)^T$ and   $g = (g_4,\ldots,g_{2N-1})$.
Consider   the  Vandermonde like  matrix
\begin{align*}
	S  
	=\begin{bmatrix}
		H_0(v_1) &&\cdots&&H_0(v_{2N})\\
		H_1(v_1)&&\cdots&&H_1(v_{2N})\\
		&&	\vdots&&\\
		H_{2N-1}(v_1)&&\cdots&&H_{2N-1}(v_{2N})\\
	\end{bmatrix}	
\end{align*}
with $S \in \R^{2N \times 2N}$ 
transforming  the  variables $f=(f_1, \ldots f_{2N})^T$ into the moments $ G= S f. $ 

Using the recursion relation of the Hermite polynomials and discrete orthogonality,
the kinetic equation is rewritten in moment variables as
\begin{align}\label{eq:macro_4eq_hbgk}
	\begin{aligned}
		\partial_t g_0 + \partial_x g_1 &= 0\\
		\partial_t g_1+ \sqrt{2} \partial_x (g_2 + \frac{1}{\sqrt{2}}g_0 )&= 0\\
		\partial_t g_{2} + \sqrt{2}\partial_x (\alpha_{3} g_{3} + \alpha_2 g_1) &= 0\\
		\partial_t g_{k} + \sqrt{2} \partial_x (\alpha_{k+1}g_{k+1} + \alpha_{k} g_{k-1})&= - \frac{1}{\epsilon}g_{k} , k=3, \ldots 2N-2\\
		\partial_t g_{2N-1}+ \sqrt{2}\partial_x (\alpha_{2N-1}g_{2N-2})&= - \frac{1}{\epsilon}g_{2N-1} . \end{aligned}
\end{align}

%
Renaming gives
\begin{align}
	\begin{aligned}
		\partial_t \rho+ \partial_x q&= 0\\
		\partial_t q+  \partial_x S &= 0\\
		\partial_t S+  a^2  \partial_x q  + \partial_x (2 \sqrt{3}  g_{3} ) &= 0\\
		\partial_t g_{3} + \sqrt{2}\partial_x (\alpha_{4} g_{4} + \alpha_3 g_2) &= - \frac{1}{\epsilon}g_{3} \\
		\partial_t g_{k} + \sqrt{2} \partial_x (\alpha_{k+1}g_{k+1} + \alpha_{k} g_{k-1})&= - \frac{1}{\epsilon}g_{k} , k=4, \ldots 2N-2\\
		\partial_t g_{2N-1}+ \sqrt{2}\partial_x (\alpha_{2N-1}g_{2N-2})&= - \frac{1}{\epsilon}g_{2N-1} \end{aligned}
\end{align}
Note that  for this system we obtain in  the limit $\epsilon \rightarrow 0$ directly  the acoustic system 
(\ref{eq:2.2}).

\subsection{The discrete layer problem}
\label{disclayer}

The corresponding discrete kinetic  layer equation is
\begin{align}
	\label{ODE}
	\begin{aligned}
		\partial_xq &= 0\\
		\partial_x S &= 0\\
		\partial_x  (  a^2 q + 2  \sqrt{3}g_{3} ) &= 0\\
		\sqrt{2}\partial_x (\alpha_{4} g_{4} + \alpha_3 g_2) &= -g_3\\
		\sqrt{2}\partial_x (\alpha_{5} g_{5} + \alpha_4 g_3) &= -g_4\\
		\sqrt{2} \partial_x (\alpha_{k+1}g_{k+1} + \alpha_{k} g_{k-1})&= - g_{k} , k=5, \ldots 2N-2\\
		\sqrt{2}\partial_x (\alpha_{2N-1}g_{2N-2})&= -g_{2N-1}
	\end{aligned}
\end{align}

That is $q=C $, $S=D$ and $D= 2 g_2 +\rho$.  and $a^2q + 2 \sqrt{3} g_3= a^2 C + 2 \sqrt{3} g_3$ a further constant. We set $g_3 =G$.
For constants $C,G \in \R$ and $D \in \R^+$.
This leads to 
\begin{align}
	\begin{aligned}
		\sqrt{2}\partial_x (\alpha_{4} g_{4} + \alpha_3 g_2) &= -G\\
		\sqrt{2}\partial_x (\alpha_{5} g_{5} )&= -g_4\\
		\sqrt{2}\partial_x (\alpha_{k+1}g_{k+1} + \alpha_{k} g_{k-1})&= - g_{k} , k=5, \ldots 2N-2\\
		\sqrt{2}\partial_x (\alpha_{2N-1}g_{2N-2})&= - g_{2N-1} \end{aligned}
\end{align}

A bounded solution is only obtained for $G=0$ and then $2 g_4 +\sqrt{3}g_2 =H $.
With $\rho= D-2 g_2$ rewriting gives 
$$
\rho = B + \frac{4}{\sqrt{3}} g_4
$$
with the new constant
$B= D-H\frac{2}{\sqrt{3}}$.
Collecting we have
\begin{align*}
	g_0 &= \frac{1}{\sqrt{2}} \rho = \frac{B}{\sqrt{2} }+ \frac{2 \sqrt{2}}{\sqrt{3}} g_4\\
	g_1 &= \frac{1}{\sqrt{2}} C \\
	g_2 &= \frac{1}{\sqrt{3}}  (H-  2  g_4)   =  \frac{1}{\sqrt{3}}  ( \frac{\sqrt{3}}{2} (D-B)-  2  g_4)  
	= \frac{D-B}{2} -\frac{2}{\sqrt{3}} g_4 \\
	g_3 &=0
\end{align*}
and the system for $g=(g_4, \ldots g_{2N-1})$  given in matrix form by 
$$
\partial_x g = - A^{-1} g 
$$
with the symmetric tridiagonal matrix $A \in \R^{2(N-2) \times 2(N-2)}$ given by 
\begin{align*}
	\begin{bmatrix}
		0 & \alpha_5& 0&0&\cdots &0\\ \alpha_5 &0&\alpha_6 & 0&\cdots &0 \\
		&\ddots &\ddots&\ddots &\\
		&&\ddots &\ddots&\ddots &\\
		0& \cdots &0&\alpha_{2N-2}&0 &\alpha_{2N-1}\\
		0& \cdots &0&0 &\alpha_{2N-1}&0
	\end{bmatrix} 
\end{align*}
The fix point of the linear ODE  system  (\ref{ODE}) is given by
$g=0$ and $\rho =B$ and $q=C$ and $S=D$. 

\begin{lemma} 
	\label{lemma1}
	$A$ is strictly hyperbolic, that means it is diagonalizable with real and distinct eigenvalues. Moreover,  $N-2$ eigenvalues of $A$ are strictly positive. The remaining $N-2$ eigenvalues have the corresponding negative values.
	We denote the eigenvectors associated to positive eigenvalues by $r_i, i=1,\ldots, N-2$ 
	and the matrix of those eigenvectors as  $$R_2^+= \begin{pmatrix} r_1,\ldots,r_{N-2}\end{pmatrix} \in \R^{2(N-2)\times (N-2)}.$$
	

\end{lemma}

\section{Solution of the coupling problem}
\label{couplingconditions}
In this section we discuss the solution of the kinetic coupling problem for the discrete velocity model.

First, we observe, that we obtain a bounded solution of the discrete kinetic half space problem if the initial values 
$(g_4(0),\cdots , g_{2N-1}(0))$ of the discrete half-space problem are located  in  the stable manifold of the equations spanned by the eigenvectors $r_1, \cdots  r_{N-2} \in \R^{2N-4}$ associated to the positive eigenvalues of $A$. That means $g$ has to fulfil

\begin{align*}
	g(0) =\gamma_1 r_1 + \cdots + \gamma_{N-2} r_{N-2} = R_2^+ \gamma
\end{align*}
for some values $\gamma_1, \ldots, \gamma_{N-2}$ in $\R$.
Then, one obtains directly with the above considerations
\begin{align}
\label{f0}
	G(0,t) = (u,g)^T (0,t) =  T (D,C,B,\gamma)^T
\end{align}
where 
%
%
$T \in \R^{(2N)  \times (N+1)}$ given by 
\begin{align*}
	T=	\begin{bmatrix}
		T_{11}&T_{12}\\
		0 & R_2^+
	\end{bmatrix}		
\end{align*}
with
\begin{align*}
	T_{11}=	\begin{bmatrix}
	 0&0&\frac{1}{\sqrt{2}} \\
		0&\frac{1}{\sqrt{2}} &0\\
		\frac{1}{2}&0&-\frac{1}{2}\\
		0&0&0
	\end{bmatrix}		,
	T_{12} = \begin{bmatrix}
	\frac{2 \sqrt{2}}{\sqrt{3}} 	e_1^T R_2^+ \\
	0\\
	 -\frac{2 }{\sqrt{3}} e_1^T R_2^+ \\
	 0 
	\end{bmatrix}	,
\end{align*}
where $e_1=(1,0,\ldots,0)^T$ is the unit vector in $\R^{2(N-2)}$.
$T$ is obviously  full  rank. Finally, one obtains 
$f (0,t)=(f_1, \ldots f_{2N})^T (0,t)$ via  the transformation $f (0,t)= S^{-1} G(0,t) $.
Using this expression for $f(0,t)$ in the kinetic coupling conditions
\begin{align}
	f^i(0,v) = \sum_{j=1}^n \beta_{ij}f^j (0,-v), v >0\ , i=1, \cdots n
\end{align}
gives $n N$ equations for $n(N+1)$ unknowns.

As discussed previously in Section \ref{kinlayer} the remaining $n$ equations are obtained from the outgoing characteristic 
quantities  
%
%
%
%
\begin{align}
	\label{cchar}
D^i - a C^i = S^i_0(0) - a q_0^i(0) , i=1, \cdots ,n.
\end{align}
This is  the final set of $n(N+1)$ equations for   the quantities $C^i, D^i,B^i, \gamma^i$ and in particular one obtains the full solution of the coupled half-space problem and 
$$
S_0^i(0) + a q_0^i(0)= D^i+a C^i
$$
i.e.  the required values for the ingoing characteristics of the  macroscopic equations at the nodes.

For  a more detailed consideration we restrict ourselves to the  case of symmetric coupling conditions.

\subsection{Symmetric coupling conditions}
In case of fully symmetric coupling conditions with $\beta_{ij} = \frac{1}{n-1}, i \neq 0$  and $\beta_{ii} =0$ 
the complexity can be strongly reduced. 
Note first that the coupling conditions 
\begin{align*}
	f^i(0,v) = \frac{1}{n-1}\sum_{j=1, j\ne i }^n f^j (0,-v), v >0\ , i=1, \cdots n
\end{align*}
give for $v>0$ and $i\ne j$
\begin{align*}
	(n-1) f^i(0,v) &= \sum_{l=1, l\ne i }^n f^l (0,-v)\\
	&=\sum_{l=1, l\ne j }^n f^l (0,-v)+f^j(0,-v) - f^i(0,-v)\\
	&= (n-1) f^j (0,v)+f^j(0,-v) - f^i(0,-v).
\end{align*}
Thus,
$$
(n-1) f (0,v)+f(0,-v) 
$$
is an invariant at the nodes and we obtain  for the discretized equations $N$ invariants  at the nodes depending on $B,C,D$ and $\gamma_1, \cdots, \gamma_{N-2}$.
\begin{align*}
	Z_1&=(n-1)f_{N+1}(0)+f_N(0)\\
	Z_i&=(n-1)f_{2N-k+1}(0)+f_k(0), k=2, \cdots N-1\\
	Z_N&=(n-1)f_{2N}(0)+f_1(0)
\end{align*}

%
%

We have with $Z=(Z_1,Z_2,\cdots, Z_N)$
\begin{align*}
	Z= 	R S^{-1} T	
	\begin{bmatrix}
		D\\C\\B\\\gamma
	\end{bmatrix}				
\end{align*}
with $R \in \R^{N\times 2N}$
\begin{align*}
	R=	\begin{bmatrix}
		0 & \cdots &0 &  0&1&n-1&0 &0&\cdots&0\\
		0 & \cdots &0&1&0&0&n-1 &0&\cdots &0\\
		&&\vdots &&\\
		1 & 0&0&&\cdots &\cdots &&0&0 &n-1
	\end{bmatrix}		
\end{align*}
Assuming that   $R S^{-1} T$ is full rank a QR decomposition gives a
transformation 
of the matrix $RS^{-1} T$ to the form
 
 \begin{align*}
 	\begin{bmatrix}
 		1 & \delta_1 & 0 &0&0&\cdots& 0 \\
 		0&1 & \frac{1}{\delta_2} & 0 &0&\cdots& 0 \\
 		0& 0&1 &\tilde \delta_1&0&\cdots&0\\
 		&&\vdots&&\\
 		0&\cdots &&\cdots&0&1&\tilde \delta_{N-2}
 	\end{bmatrix}		
 \end{align*}
 
 
In particular, this  gives the invariants
\begin{align*}
	D +  \delta_1 C
\end{align*}
and
\begin{align*}
B +  \delta_2 C
\end{align*}
or in terms of $S$ and $q$ and $\rho$
the invariance of 
\begin{align*}
	S_\infty + \delta_1 q_\infty
\end{align*}
and
\begin{align*}
\rho_\infty + \delta_2 q_\infty
\end{align*}
This gives  as discussed above $2(n-1)$  equations at each  node, i.e.
the conditions
\begin{align*}
	\label{cinv}
	D^i +  \delta_1 C^i =D^j +  \delta_1 C^j  , \; 
	B^i +  \delta_2 C^i=B^j +  \delta_2 C^j.
\end{align*}

Together with the equality of fluxes 
\begin{align}
	\label{cflux}
\sum_{i=1}^n C^i =0
\end{align}
we have therefore $2n-1$
coupling conditions.
Additionally we obtain  $n$ more conditions from conditions (\ref{outcond}) as before.
This gives altogether again $3n-1$ equations for the $3n$ unknown quantities $C^i$ and $D^i$ and $B^i$  at each node.
the final condition (\ref{balancerho}) is obtained from the consideration of the viscous layer.

\begin{remark}
We note that  already from the invariance of $S_{\infty}+\delta_1q_{\infty}$, the  balance of fluxes $\sum_j q_{\infty}^{j}=0 $
and the $n$ conditions at infinity, we have a sufficient number of conditions for the coupling of the macroscopic equations: 
this fixes the values  of $ q_{\infty}^i$ and $S_{\infty}^i$ at the node, which is sufficient to couple the macroscopic equations due to the $0$-characteristic
not requiring a boundary condition. 
\end{remark}

\subsection{Full kinetic solution at the node}
\label{nearnode}

To obtain the full  kinetic solution at the node in the limit $\epsilon \rightarrow 0 $  we have to 
determine the layer solution at $x=0$. That means according to (\ref{f0}) we have to determine the values of $\gamma_1, \ldots, \gamma_{N-1}$. That gives finally all moments of the distribution function at the node.
In case of fully symmetric coupling conditions we can simplify the procedure.
Note that in this case   the transformation of the matrix $Z$ gives  the additional $N-2$ invariants
\begin{align}
	\label{invnext}
	\begin{aligned}
		C+ \tilde \delta_1 \gamma_1, \\
		\gamma_{k-1} + \tilde \delta_{k} \gamma_{k}, \; k=2, \ldots N-2 .
	\end{aligned}
\end{align}
Moreover, we obtain directly from the coupling conditions for the odd moments
\begin{align*}
	\sum_{i=1}^n g^i_{2k+1}(x=0) =0, k= 2, \ldots, N-1,
\end{align*}
which leads to 
\begin{align}
	\label{invnext2}
	\sum_{i=1}^n  e_{2k}^T R_2^+ \gamma^i =0, k= 2, \ldots, N-1 .
\end{align}
(\ref{invnext}) and (\ref{invnext2})  give the required $(N-2)(n-1)+N-2= (N-2)n$ conditions additionally to the $3n$ conditions from above and therefore 
$C^i,D^i, B^i,\gamma^i$ and thus all moments $\rho^i,q^i, S^i, g^i$ at $x=0$.
In particular,  $\rho^i(x=0)$ is given by 
$$\rho^i (x=0)= B+\frac{4}{\sqrt{3}} g^i_4= B+\frac{4}{\sqrt{3}} e_1^T R_2^+   \gamma^i$$.

The  distribution functions on edge $i$  at the node
 $f^i=f^i (x=0, v) $ for $v \in \R$ is finally given by the Hermite expansion

\begin{align*}
	f(v) = H_0(\frac{v}{\sqrt{2}})   \sum_{k=0}^{2N-1} g_k H_k (\frac{v}{\sqrt{2}})
\end{align*}
where
\begin{align*}
	g_0= \frac{\rho}{\sqrt{2}}, g_1 =\frac{C}{\sqrt{2}}, g_2= \frac{1}{2} (D-\rho), g_3 =0
\end{align*}
and  $g =(g_4, \ldots,g_{2N-1})^T$ as before as
\begin{align*}
	g =  R_2^+ \gamma.
\end{align*}

\subsection{Summary of coupling conditions for the macroscopic quantities}
\label{sum}
The above   3n coupling conditions (\ref{cinv}) (\ref{cflux}) (\ref{cchar}) and (\ref{balancerho}) can be written as a  linear system  $	\mathcal{A} m = b $ for
\begin{align*} 
 m&=(S^1_\infty, \cdots S^n_\infty, q^1_\infty \cdots q^n_\infty,\rho^1_\infty, \cdots \rho^n_\infty) \\
 &=  (D^1, \cdots,D^n, C^1, \cdots, C^n, B^1, \cdots, B^n)
\end{align*}
 given by  the block matrix 

\begin{align*} 
	\mathcal{A}= \begin{bmatrix}
		A_{11}&A_{12}&0\\
		A_{21}&A_{22}&0\\
		A_{31}&A_{32}&A_{33}
	\end{bmatrix} \in \R^{3n\times 3n}
\end{align*}
with $A_{ij} \in \R^{n \times n}$ defined by 
$$
A_{11} = A, A_{12} = B+\delta_1 A, A_{31} = B,A_{32} = \delta_2 A, A_{33} =-3B+A
$$
and
$$
A_{21} = I,  A_{22} =- a I,
$$
where we have defined 
\begin{align*} 
A= \begin{bmatrix}
0&0&0& \cdots &0&0\\
1&-1&0& 0&\cdots &0\\
0&1&-1 &0&\cdots &0\\&\vdots & \vdots&\\0&\cdots &0&0&1&-1
\end{bmatrix} \in \R^{n \times n}
\end{align*}
and
\begin{align*} 
B= \begin{bmatrix}
1&  1&\cdots &1&1\\
0&0&0&\cdots &0\\& \vdots&\\0&\cdots &0&0&0
\end{bmatrix} \in \R^{n \times n}
\end{align*}
and 
\begin{align*} 
b= \begin{bmatrix}
	0^n   &\alpha &\beta&0^{n-1}
\end{bmatrix}^T \in \R^{3 n}.
\end{align*}
Here, we have used $\alpha = (\alpha^{(1)}, \cdots, \alpha^{(n)} )\in \R^n$ with $ \alpha_i = S^{(i)}_0- a  q^{(i)}_0$ and $  \beta = \sum_{i=1}^n [S^i_0 -3 \rho^i_0] \in \R$. Moreover,  $0^n \in \R^n $ is the $0$-vector.


$\mathcal{A}$ is invertible, that means the coupling problem is uniquely solvable, if the 
determinant of $\mathcal{A}$ given by 
\begin{align*}
	\det(A) &= \det(A_{33})\det(A_{22})  \det(A_{11} - A_{12} A_{22}^{-1} A_{21})  \\&= 
	\det(A_{33})\ \det(A_{11} +a  A_{12} )\\&=-3 a (-1)^{n+1} n    (-1)^{n+1} n (1+a \delta_1)^{n-1}   \\
	&= -3an^2   (1+a \delta_1)^{n-1} 
\end{align*}
is not equal to zero, i.e. for all $\delta_1 \neq - \frac{1}{a}$ and, in particular for all non-negative values of $\delta_1$
and arbitrary values of $\delta_2$.

\subsection{Simple approximations for  coupling conditions}
\label{halfspacemarshakbgk}

There is an extensive literature on approximate  methods to find the asymptotic states of half-space problems \cite{M47,GK95,LF2,ST,SO,LLS16}. We use two of them adapted to the network situation, see \cite{ABEK22},  for numerical comparison in Section
\ref{numericalbgk}.
	A straightforward approximation of the asymptotic states depending on the ingoing distribution is obtained by equalizing the first two half moments  at $x=0$ with  those at $x=\infty$, i.e. 
	$$
	\int_{0}^\infty \left(
	\begin{matrix}
 v\\
	v^2	
\end{matrix}
\right)  f (0,v) dv=	
\int_{0}^\infty \left(
\begin{matrix}
	v\\
	v^2	
\end{matrix} \right)
	\left(\rho+vq+\frac{1}{2}(v^2-1)(S-\rho)  \right)M(v) dv 
$$	
and assuming  that the outgoing distribution is  approximated by
$$
f(0,v)  =  \left( \rho_\infty + v q_\infty +\frac{1}{2}(v^2-1)(S_\infty-\rho_\infty) \right) M(v) , v <0\ .
$$


Using the coupling conditions 

\begin{align*}
f^i(0,v) = \sum_{j=1}^n \beta_{ij}f^j (0,-v), v >0.
\end{align*}
gives
\begin{align*}	
\int_{0}^\infty \left(
\begin{matrix}
v\\
v^2	
\end{matrix} \right)
\left(\rho+vq+\frac{1}{2}(v^2-1)(S-\rho)  \right)M(v) dv \\
=
\sum_{j=1}^n \beta_{ij}    \int_{0}^\infty \left(
\begin{matrix}
v\\
v^2	
\end{matrix} \right)\left( \rho_\infty - v q_\infty +\frac{1}{2}(v^2-1)(S_\infty-\rho_\infty) \right) M(v)   dv
\end{align*}	
or

%
%
\begin{align*}
\rho^i_\infty +  \sqrt{2\pi}q^i_\infty +S_\infty
= \sum_{j} \beta_{ij} \left( \rho^j_\infty -   \sqrt{2\pi}q^j_\infty   + S_\infty\right)
\end{align*}
and
\begin{align*}
4q^i_\infty +  \sqrt{2\pi}S_\infty
= \sum_{j} \beta_{ij} \left(
-   4q^j_\infty   +  \sqrt{2\pi}S_\infty\right)
\end{align*}

In the special case of a  uniform node, i.e. $c_{ij} = \frac{1}{n-1}$ for $i \neq j$ and $c_{ij}=0$ for $i =j$  we have

%
%
\begin{align*}
	S^i_\infty +   \frac{4(n-2)}{n\sqrt{2\pi}}  q^i_\infty 
	=
	\rho^j_\infty +   \frac{4(n-2)}{n\sqrt{2\pi}} q^j_\infty  
\end{align*}
and  the invariant
\begin{align}
	S_\infty +   \frac{4(n-2)}{n\sqrt{2\pi}}  q_\infty \ .
\end{align}
or
\begin{align}
	\delta_1 = \frac{4(n-2)}{n\sqrt{2\pi}}
\end{align}

Moreover we obtain the invariant
\begin{align}
\rho_\infty +   \frac{n-2}{n} \frac{2(\pi-2) }{\sqrt{2\pi}}q_\infty \ .
\end{align}
or
\begin{align}
\delta_2 = \frac{n-2}{n} \frac{2(\pi-2) }{\sqrt{2 \pi}}
\end{align}
For example for $n=3$, we obtain for  the  factor 
$\delta_1 = \frac{4}{3\sqrt{2\pi}}\sim 0.5320$ and $\delta_2 = \frac{4}{3}\frac{2(\pi-2) }{\sqrt{2 \pi}}\sim0.3033$. 
$n \rightarrow \infty$ gives $\delta_1 = \frac{4}{\sqrt{2\pi}} \sim 1.5958$ and $\delta_2 = \frac{2(\pi-2) }{\sqrt{2 \pi}}\sim 0.9109$.



The  half moment approximation from \cite{ABEK22} gives the values
$\delta_1\approx0.5301$ and  $\delta_2\approx0.3402$. 
For the case $n=\infty$, this yields $\delta_1\approx1.5833$ and  $\delta_2\approx0.9975$.

%
%

\subsection{Numerical results for the coupling coefficients}
\label{ccnum}
We restrict ourselves to fully symmetric coupling conditions.
From a numerical point of view the computation of $\delta_1,\delta_2$ is independent from the solution of the network problem.
It requires in particular the knowledge of the positive eigenvalues $\lambda_i, i= 1 , \ldots , N-1$ of the matrix $A$.
Moreover, an inversion of the Vandermonde like matrix $S$  is needed and one Gaussian elimination
of $RS^{-1} T$.
The inversion of $S$ is not ill conditioned, as long as the Gauß-Hermite  points are used, see, e.g. \cite{Gautschi}.

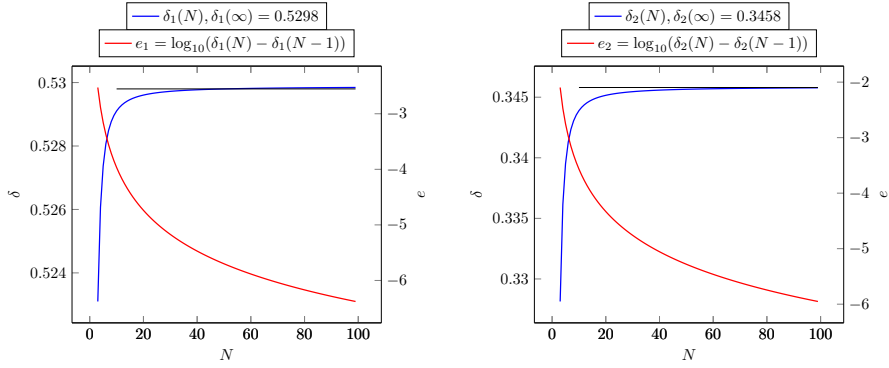
\begin{figure}[h!]
\center
	\externaltikz{netyyy}{	
		\begin{tikzpicture}[scale=0.6]
		\begin{axis}[ ylabel = $\delta$, ylabel near ticks, 
		xlabel =  $N$, xlabel near ticks,
		legend style = {at={(0.5,1)},xshift=0.0cm,yshift=0.8cm,anchor=south},
		legend columns= 1,
		yticklabel style={/pgf/number format/.cd,fixed,precision=3},
		axis y line*=left,
		]
		\addplot[color = blue, thick] file{data/datapaper/delta1_bgk3.txt};
		\addplot[black] coordinates {(10,0.5298)  (99,0.5298)};	
				\addlegendentry{$\delta_1(N), \delta_1(\infty)= 0.5298$}						
		\end{axis}
			\begin{axis}[ylabel =$e$,ylabel near ticks, 
		legend style = {at={(0.5,1)},xshift=0.0cm,yshift=0.2cm,anchor=south},
		legend columns= 1,
		axis y line*=right,
		]
		
		\addplot[color = red,thick] file{data/datapaper/error1_bgk3.txt};
		\addlegendentry{$e_1=\log_{10}(\delta_1(N) - \delta_1(N-1))$}
	
		\end{axis}
		\end{tikzpicture}
	}
\externaltikz{netyyybbb}{	
	\begin{tikzpicture}[scale=0.6]
\begin{axis}[ ylabel = $\delta$, ylabel near ticks, 
xlabel =  $N$, xlabel near ticks,
legend style = {at={(0.5,1)},xshift=0.0cm,yshift=0.8cm,anchor=south},
legend columns= 1,
yticklabel style={/pgf/number format/.cd,fixed,precision=3},
axis y line*=left,
]
\addplot[color = blue, thick] file{data/datapaper/delta2_bgk3.txt};
\addplot[black] coordinates {(10,0.3458)  (99,0.3458)};	
\addlegendentry{$\delta_2(N), \delta_2(\infty)= 0.3458$}						
\end{axis}
\begin{axis}[ylabel =$e$,ylabel near ticks, 
legend style = {at={(0.5,1)},xshift=0.0cm,yshift=0.2cm,anchor=south},
legend columns= 1,
axis y line*=right,
]

\addplot[color = red,thick] file{data/datapaper/error2_bgk3.txt};
\addlegendentry{$e_2=\log_{10}(\delta_2(N) - \delta_2(N-1))$}

\end{axis}
\end{tikzpicture}
}
	\caption{Coefficients $\delta_1  $ and $\delta_2$  depending on $N$ for $n=3$. Associated error depending on $N$.}
	\label{fig4}
\end{figure}

\begin{figure}[h!]
	\center
		\externaltikz{funcr11mmm1}{
	\begin{tikzpicture}[scale=0.6]	
	\begin{axis}[ ylabel = $\delta$, ylabel near ticks, 
	xlabel =  $N$, xlabel near ticks,
	legend style = {at={(0.5,1)},xshift=0.0cm,yshift=0.8cm,anchor=south},
	legend columns= 1,
	yticklabel style={/pgf/number format/.cd,fixed,precision=3},
	axis y line*=left,
	]
	\addplot[color = blue, thick] file{data/datapaper/delta1_bgk.txt};
	\addplot[black] coordinates {(10,1.5826)  (99,1.5826)};	
	\addlegendentry{$\delta_1(N), \delta_1(\infty)= 1.5826$}						
	\end{axis}
	\begin{axis}[ylabel =$e$,ylabel near ticks, 
	legend style = {at={(0.5,1)},xshift=0.0cm,yshift=0.2cm,anchor=south},
	legend columns= 1,
	axis y line*=right,
	]
	
	\addplot[color = red,thick] file{data/datapaper/error1_bgk.txt};
	\addlegendentry{$e_1=\log_{10}(\delta_1(N) - \delta_1(N-1))$}
	
	\end{axis}
	\end{tikzpicture}
}
\externaltikz{funcr1vvcc1}{
	\begin{tikzpicture}[scale=0.6]

	\begin{axis}[ ylabel = $\delta$, ylabel near ticks, 
	xlabel =  $N$, xlabel near ticks,
	legend style = {at={(0.5,1)},xshift=0.0cm,yshift=0.8cm,anchor=south},
	legend columns= 1,
	yticklabel style={/pgf/number format/.cd,fixed,precision=3},
	axis y line*=left,
	]
	\addplot[color = blue, thick] file{data/datapaper/delta2_bgk.txt};
	\addplot[black] coordinates {(10,1.0079)  (99,1.0079)};	
	\addlegendentry{$\delta_2(N), \delta_2(\infty)= 1.0079$}						
	\end{axis}
	\begin{axis}[ylabel =$e$,ylabel near ticks, 
	legend style = {at={(0.5,1)},xshift=0.0cm,yshift=0.2cm,anchor=south},
	legend columns= 1,
	axis y line*=right,
	]
	
	\addplot[color = red,thick] file{data/datapaper/error2_bgk.txt};
	\addlegendentry{$e_2=\log_{10}(\delta_2(N) - \delta_2(N-1))$}
	
	\end{axis}
	\end{tikzpicture}
}
	\caption{Coefficients $\delta_1  $ and $\delta_2$  depending on $N$ for $n=\infty$. Associated error depending on $N$.}
	\label{fig5}
\end{figure}
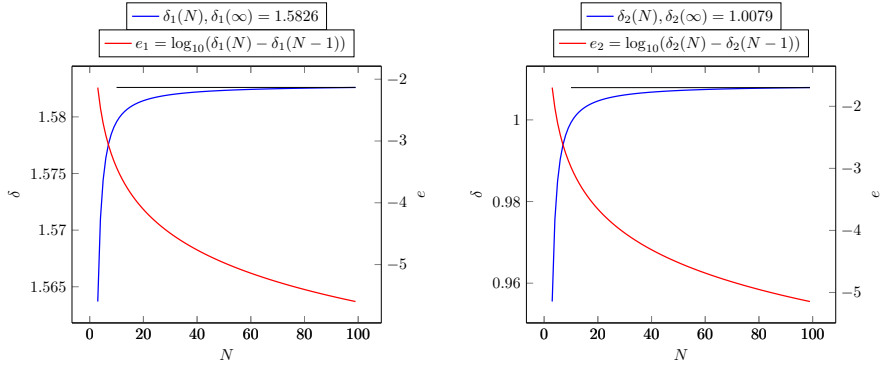

Figure \ref{fig4} shows the resulting values of $\delta_1,\delta_2$ for different values of $N$ together with
the associated numerical error for $n=3$ edges at the node.
Figure \ref{fig5} shows the same, but for $n=\infty$ edges at the node.

\section{Numerical results for the network}
\label{numericalbgk}

To illustrate the above results,  we consider the case of  a single node  with 3 edges.

As initial conditions for the kinetic and macroscopic equations we choose the 
 densities $\rho_0= (\rho^1_0,\rho^2_0,\rho^3_0) = (1,1-\bar \rho_0,1+\bar \rho_0)$,  the fluxes $q_0 =( 0, \bar q_0  , -\bar q_0 )$ and the second moment $S_0= (1,1-\bar S_0, 1+\bar S_0)$ and   the associated linearized Maxwellians for the kinetic solution.
 These initial  data are also prescribed at the  outer boundaries.
 
 The coupling conditions from Section  \ref{sum} yield the following solutions.
 \begin{align*}
 q_\infty&=(0,\bar q_\infty , -\bar q_\infty)\\
 S_\infty &= (1,1-\bar S_\infty, 1+ \bar S_\infty)\\
 \rho_\infty&=(1, 1-\bar \rho_\infty, 1+\bar \rho_\infty)
 \end{align*} 
 with $$\bar q_\infty = \frac{\bar S_0 +a \bar q_0}{\delta_1 +a } , \; \bar S_\infty = \delta_1 \bar q_\infty   , \; \tilde \rho = \delta_2  \bar q_\infty.$$

Note that in the variable $\rho$ kinetic and viscous  layer may  appear  additionally to the travelling wave solution, whereas in $q$ and $S $ only internal travelling wave solutions exist. The later quantities are constant throughout kinetic and viscous layer.
See the considerations in section \ref{nnode}.
In the following we discuss different physical situations. 

Case 1 is a situation without internal travelling wave solution and without viscous layer.
For this, we choose initial conditions with $\bar q_0=1$, $\bar \rho_0 = \delta_2 \bar q_0$ and $\bar S_0 = \delta_1 \bar q_0$, such that the initial conditions fulfil the coupling conditions, since 
$\bar q_\infty =\bar q_0=1$, $\bar S_\infty=  \bar S_0=\delta_1 $ and $\bar \rho_\infty = \bar \rho_0=\delta_2  $. In this case only a kinetic layer in $\rho$ appears at the node. This situation is shown in Figure \ref{fig_n1}, where on the right a zoom to the region near the node is shown.

Case 2 is a case with merged kinetic and viscous layer in $\rho$ and no travelling wave in $q $ and $S$.
For this we choose $\bar q_0 =1 $, $\bar S_0= \delta_1 \bar q_0$, but $\bar \rho_0 = 2 \delta_2 $. 
We obtain $\bar q_\infty =\bar q_0=1$, $\bar S_\infty =  \bar S_0 =\delta_1 $ and $\bar \rho_\infty = \delta_2 \bar q_0  =  \delta_2 \neq \bar \rho_0 $.
This   avoids travelling wave solutions, but leads to a viscous layer near the node
mixing with the kinetic layer.
This situation is shown in Figure \ref{fig_n2}, where again on the right a zoom  to the region near the node is shown.

Since there is no travelling wave in the zero characteristic, we have according to the considerations in Section \ref{nnode}
a left travelling wave state in $\rho$ given by   $ \rho_0+ \frac{1}{3} (  S_\infty  - S_0 )  $. Thus, if  $\rho_\infty =  \rho_0+ \frac{1}{3} (  S_\infty  - S_0 )  $ or equivalently $\bar \rho_\infty = \bar \rho_0+ \frac{1}{3} (  \bar S_\infty  -  \bar S_0 )  $, then,  there is no viscous layer.
Accordingly, for test case 3, we choose  $\bar q_0 =1 $, $\bar S_0 = 2 \delta_1 \bar q_0$ and  $$\bar \rho_0 = \delta_2 \frac{ 2 \delta_1 +a }{\delta_1 +a } - \frac{1}{3} (  \delta_1 \frac{ 2 \delta_1 +a }{\delta_1 +a }  -  2 \delta_1 )= \frac{ 2 \delta_1 +a }{\delta_1 +a }(\delta_2 -\frac{1}{3}  \delta_1) + \frac{2}{3}\delta_1.$$
We obtain $\bar q_\infty =\frac{ 2 \delta_1 +a }{\delta_1 +a } $, $\bar S_\infty =  \delta_1 \bar q_\infty $ and $\bar \rho_\infty = \delta_2 \bar q_\infty  =\delta_2 \frac{ 2 \delta_1 +a }{\delta_1 +a } = \bar \rho_0+ \frac{1}{3} (\bar  S_\infty  -  \bar S_0 ) $.
This yields  travelling waves in $q$ and $S$ and a  kinetic layer, but no viscous layer  in $\rho$ as can be seen in Figure \ref{fig_n3}.

Finally, for test case 4 we choose $\bar q_0 =1 $, $\bar S_0 = 2 \delta_1 \bar q_0$ and  $\bar \rho_0 = \delta_2 \frac{2 \delta_1 +a}{\delta_1 +a } $. 
We obtain $\bar q_\infty =\frac{ 2 \delta_1 +a }{\delta_1 +a } $, $\bar S_\infty =  \delta_1 \bar q_\infty $ and $\bar \rho_\infty = \delta_2 \bar q_\infty = \bar \rho_0  $.
This gives  travelling waves in $q$ and $S$ and in $\rho$  a merged kinetic and  viscous layer
as can be seen in Figure \ref{fig_n4}.

Finally Figure \ref{fig_nf} shows the full distribution functions at $x=0$ for test case 2. In particular, one observes the discontinuity at $v=0$
and the Gibbs oscillations for the spectral method.

\begin{figure}		
	\center									
	\externaltikz{nebb2}{
		\begin{tikzpicture}[scale=0.6]
			\begin{axis}[ylabel = $\rho$,xlabel = $x$,
				legend style = {at={(0.5,1)},xshift=0.2cm,yshift=-0.0cm,anchor=south},
				legend columns= 2,
				xmin = 0.0, xmax = 0.1,	
				]

					del1=0.5298;
				del2=0.3458;
				
				rhoinf=0.6542;

				\addplot[color = brown,thick] file{network_matlab/testcase1/rho_kinetic_1.txt};
				\addlegendentry{kinetic $1$ };

					\addplot[mark=none, black, thick, domain = 0:0.01] {0.6542};
				\addlegendentry{spectral $\rho^2_\infty$};
				
				\addplot[color = blue!,thick] file{network_matlab/testcase1/rho_kinetic_2.txt};
				\addlegendentry{kinetic $2$};

				\addplot[mark=none, thick,red,domain = 0:0.01]{0.7245};
				\addlegendentry{spectral $\rho^2(0)$};

					\addplot[color = green,thick] file{network_matlab/testcase1/rho_kinetic_3.txt};
				\addlegendentry{kinetic $3$};

			\end{axis}
		\end{tikzpicture}
	}									
			\externaltikz{nxxkxff}{
			\begin{tikzpicture}[scale=0.6]
			\begin{axis}[ylabel = $\rho$,xlabel = $x$,
				legend style = {at={(0.5,1)},xshift=0.2cm,yshift=-0.0cm,anchor=south},
				legend columns= 1,
				xmin = 0.0, xmax = 0.008,
				ymin = 0.6, ymax=0.8,	
				]
				
				\addplot[color = blue!,thick] file{network_matlab/testcase1/rho_kinetic_2.txt};
				\addlegendentry{kinetic $2$};
				
				\addplot[mark=none, black, thick, domain = 0:0.01] {0.6542};
				\addlegendentry{spectral $\rho_\infty^2$};

				\addplot[mark=none, thick,red,domain = 0:0.0005]{0.7245};
	\addlegendentry{spectral $\rho^2(0)$};
			\end{axis}
		\end{tikzpicture}
	}	
\hspace{1cm}
	\caption{Testcase 1: $\rho$ for all edges, kinetic solution for   $\epsilon= 5 \cdot 10^{-4}$ at time $t=0.1$ (left). Zoom to solution on edge 2 (right).}
	\label{fig_n1}
\end{figure}

\begin{figure}		
	\center									
	\externaltikz{nggg2}{
		\begin{tikzpicture}[scale=0.6]
			\begin{axis}[ylabel = $\rho$,xlabel = $x$,
				legend style = {at={(0.5,1)},xshift=0.2cm,yshift=-0.0cm,anchor=south},
				legend columns= 2,
				xmin = 0.0, xmax = 0.1,	
				]
				
					del1=0.5298;
				del2=0.3458;
				
				rhoinf=0.6542;
				\addplot[color = brown,thick] file{network_matlab/testcase2c/rho_kinetic_1.txt};
				\addlegendentry{kinetic $1$ };
				
					\addplot[mark=none, black, thick, domain = 0:0.1] {0.6542};
				\addlegendentry{spectral $\rho_\infty^2$};
				
				\addplot[color = blue!,thick] file{network_matlab/testcase2c/rho_kinetic_2.txt};
				\addlegendentry{kinetic $2$};
				
					\addplot[mark=none, thick,red,domain = 0:0.01]{0.7245};
				\addlegendentry{spectral $\rho^2(0)$};

				\addplot[color = green,thick] file{network_matlab/testcase2c/rho_kinetic_3.txt};
				\addlegendentry{kinetic $3$};

			\end{axis}
		\end{tikzpicture}
	}									
	\externaltikz{nfiigf}{
		\begin{tikzpicture}[scale=0.6]
			\begin{axis}[ylabel = $\rho$,xlabel = $x$,
				legend style = {at={(0.5,1)},xshift=0.2cm,yshift=-0.0cm,anchor=south},
				legend columns= 2,
				xmin = 0.0, xmax = 0.003,
				ymin = 0.5, ymax=0.9,	
				]
				
%
%
					\addplot[mark=none, black, thick, domain = 0:0.04] {0.6542};
				\addlegendentry{ spectral $\rho_\infty^2$};
				
				
					\addplot[mark=none, thick,red,domain = 0:0.0003]{0.7245};
				\addlegendentry{ spectral $\rho^2(0)$};
				
					\addplot[color = blue!,thick] file{network_matlab/testcase2e/rho_kinetic_2.txt};
				\addlegendentry{kinetic $\epsilon=  10^{-4}$};

			\end{axis}
		\end{tikzpicture}
	}	
%
%
%
%
%
%
	\hspace{1cm}
	\caption{Testcase2: $\rho$ for all edges, kinetic solution for   $\epsilon=   \cdot 10^{-4}$ at time $t=0.1$ (left). Zoom to solution on edge 2 (right).}
	\label{fig_n2}
\end{figure}

\begin{figure}		
	\center									
	\externaltikz{neccx2}{
		\begin{tikzpicture}[scale=0.6]
		\begin{axis}[ylabel = $\rho$,xlabel = $x$,
		legend style = {at={(0.5,1)},xshift=0.2cm,yshift=-0.0cm,anchor=south},
		legend columns= 2,
		xmin = 0.0, xmax = 0.3,	
		]

		\addplot[color = brown,thick] file{network_matlab/testcase3b/rho_kinetic_1.txt};
		\addlegendentry{kinetic $1$ };
		
		\addplot[color = blue!,thick] file{network_matlab/testcase3b/rho_kinetic_2.txt};
		\addlegendentry{kinetic $2$};
		\addplot[color = green,thick] file{network_matlab/testcase3b/rho_kinetic_3.txt};
		\addlegendentry{kinetic $3$};
		
		del1=0.5298;
		del2=0.3458;
		
		rhoinf=0.6542;
		\addplot[mark=none, black, thick, domain = 0:0.1] {0.5732};
		\addlegendentry{spectral $\rho_\infty^2$};
		
		\addplot[mark=none, thick,red,domain = 0:0.01]{0.6599};
			\addlegendentry{spectral $\rho^2(0)$};
		
		\end{axis}
		\end{tikzpicture}
	}									
	\externaltikz{nexkpf}{
		\begin{tikzpicture}[scale=0.6]
		\begin{axis}[ylabel = $\rho$,xlabel = $x$,
		legend style = {at={(0.5,1)},xshift=0.2cm,yshift=-0.0cm,anchor=south},
		legend columns= 1,
		xmin = 0.0, xmax = 0.04,
		ymin = 0.2, ymax=0.8,	
		]
		
		\addplot[color = blue!,thick] file{network_matlab/testcase3b/rho_kinetic_2.txt};
		\addlegendentry{kinetic $2$};
		
		\addplot[mark=none, black, thick, domain = 0:0.04] {0.5732};
		\addlegendentry{ spectral $\rho_\infty^2$};

		\addplot[mark=none, thick,red,domain = 0:0.005]{0.6599};
		\addlegendentry{ spectral $\rho^2(0)$};
		\end{axis}
		\end{tikzpicture}
	}	
	\hspace{1cm}
	\caption{Testcase 3: $\rho$ for all edges, kinetic solution for    $\epsilon= 5 \cdot 10^{-4}$ at time $t=0.1$ (left). Zoom to solution on edge 2 (right).}
	\label{fig_n3}
\end{figure}

\begin{figure}		
	\center									
	\externaltikz{necdxxxddfc2}{
		\begin{tikzpicture}[scale=0.6]
		\begin{axis}[ylabel = $\rho$,xlabel = $x$,
		legend style = {at={(0.5,1)},xshift=0.2cm,yshift=-0.0cm,anchor=south},
		legend columns= 2,
		xmin = 0.0, xmax = 0.3,	
		]

		\addplot[color = brown,thick] file{network_matlab/testcase4/rho_kinetic_1.txt};
		\addlegendentry{kinetic $1$ };
		
		\addplot[color = blue!,thick] file{network_matlab/testcase4/rho_kinetic_2.txt};
		\addlegendentry{kinetic $2$};
		\addplot[color = green,thick] file{network_matlab/testcase4/rho_kinetic_3.txt};
		\addlegendentry{kinetic $3$};
		
		del1=0.5298;
		del2=0.3458;
		
		rhoinf=0.6542;
		\addplot[mark=none, black, thick, domain = 0:0.3] {0.5732};
		\addlegendentry{spectral $\rho_\infty^2$};
		
		\addplot[mark=none, thick,red,domain = 0:0.01]{0.6599};
			\addlegendentry{spectral $\rho^(0)2$};

		\end{axis}
		\end{tikzpicture}
	}									
	\externaltikz{npxxxdfdf}{
		\begin{tikzpicture}[scale=0.6]
		\begin{axis}[ylabel = $\rho$,xlabel = $x$,
		legend style = {at={(0.5,1)},xshift=0.2cm,yshift=-0.0cm,anchor=south},
		legend columns= 1,
		xmin = 0.0, xmax = 0.04,
		ymin = 0.2, ymax=0.8,	
		]
		
		\addplot[color = blue!,thick] file{network_matlab/testcase4/rho_kinetic_2.txt};
		\addlegendentry{kinetic $2$};
		
		\addplot[mark=none, black, thick, domain = 0:0.04] {0.5732};
		\addlegendentry{ spectral $\rho_\infty^2$};

		\addplot[mark=none, thick,red,domain = 0:0.005]{0.6599};
		\addlegendentry{ spectral $\rho^2(0)$};
		\end{axis}
		\end{tikzpicture}
	}	
	\hspace{1cm}
	\caption{Testcase 4: $\rho$ for all edges, kinetic solution for    $\epsilon= 5 \cdot 10^{-4}$ at time $t=0.1$ (left). Zoom to solution on edge 2 (right).}
	\label{fig_n4}
\end{figure}

\begin{figure}											
	\externaltikz{net2}{
		\begin{tikzpicture}[scale=0.7]
			\begin{axis}[ylabel = $f$,xlabel = $v$,
				legend style = {at={(0.5,1)},xshift=0.2cm,yshift=-0.0cm,anchor=south},
				legend columns= 2,
				xmin = -5, xmax = 5,
				y tick label style={/pgf/number format/precision=4},
				]
				
					\addplot[color = brown,thick] file{data/f_hermite1.txt};
				\addlegendentry{spectral $1$}
					\addplot[color = brown,thick] file{network_matlab/testcase2b/f_kinetic_1.txt};
				\addlegendentry{kinetic $1$}
				\addplot[color = red,thick] file{data/f_hermite2.txt};
				\addlegendentry{spectral $2$}
					\addplot[color = blue,thick] file{network_matlab/testcase2b/f_kinetic_2.txt};
				\addlegendentry{kinetic $2$}
				\addplot[color = magenta,thick] file{data/f_hermite3.txt};
				\addlegendentry{spectral $3$}

				\addplot[color = green,thick] file{network_matlab/testcase2b/f_kinetic_3.txt};
				\addlegendentry{kinetic $3$}

				\addplot [thick] coordinates {(0.0,0.0) (0.0,0.6)};
				%
				%

				
			\end{axis}
		\end{tikzpicture}
	}
			\externaltikz{net1}{
			\begin{tikzpicture}[scale=0.7]
				\begin{axis}[ylabel = $f$,xlabel = $v$,
					legend style = {at={(0.5,1)},xshift=0.2cm,yshift=-0.0cm,anchor=south},
					legend columns= 1,
					xmin = -0.6, xmax = 0.6,
					ymin = 0.2, ymax=0.5,
					y tick label style={/pgf/number format/precision=4},
					]
					
						\addplot[color = red,thick] file{data/f_hermite2.txt};
					\addlegendentry{spectral $2$}
					\addplot[color = blue,thick] file{network_matlab/testcase2b/f_kinetic_2.txt};
					\addlegendentry{kinetic $2$}

					
					\addplot [thick] coordinates {(0.0,0.0) (0.0,0.6)};
					%

					
				\end{axis}
			\end{tikzpicture}
	}
	\caption{Kinetic solution for all edges at $x=0$   at time $t=0.1$ for testcase 2 computed by FD method with   $\epsilon= 5 \cdot 10^{-4}$ and $\Delta x = 10^{-4}, \Delta v = 0.005$ and by spectral method with $N=1000$. On the right a zoom to the discontinuity at $v=0$ is shown for edge 2.}
	\label{fig_nf}
\end{figure}
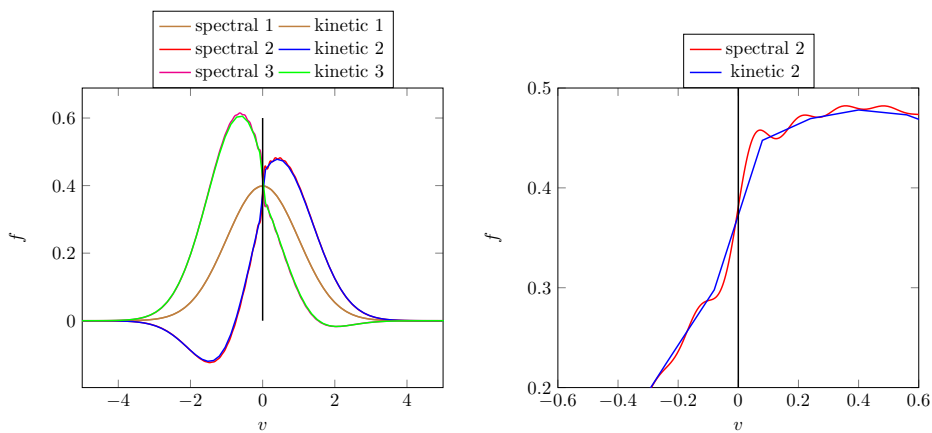
 \section{Conclusion and Outlook}
 
 This paper gives a detailed numerical analysis of the interface layers at a node of a network with dynamics described by the full  linearised kinetic BGK equations and its macroscopic limit, the linearised Euler equations.
 
 The procedure is based on an asymptotic analysis  of the situation near the nodes and the investigation of 
 the kinetic layer near the nodes and the associated coupled kinetic half-space problems. Moreover, 
 since the macroscopic equations are degenerate in one variable, the investigation requires, additionally to the discussion of the kinetic half-space problems also the investigation of related viscous layers.
 
 For the numerical solution an accurate  spectral procedure to solve the half-space problem and determine the asymptotic states  has been developed.

The validity of the  asymptotic expansion and convergence of the kinetic solution towards the macroscopic solution is proven in a following paper, see
\cite{ZK} using the  approach from \cite{WY} for hyperbolic relaxation systems.

Finally, we mention, that codes and data that allow readers to reproduce the most important numerical results in this paper, in particular the computation of $\delta_1, \delta_2 $ and the distribution function at the node,  are available at 

https://gitlab.rhrk.uni-kl.de/klar/kinetic-network-degenerate.git



\end{document}